\documentclass[11pt, a4paper]{article}
\usepackage{array}
\usepackage{multirow}
\makeatletter
\newcommand{\thickhline}{%
    \noalign {\ifnum 0=`}\fi \hrule height 1pt
    \futurelet \reserved@a \@xhline
}
\newcolumntype{"}{@{\hskip\tabcolsep\vrule width 1pt\hskip\tabcolsep}}
\makeatother

\usepackage{graphics}
\usepackage{graphicx}
\usepackage{epsfig}
\usepackage{subfigure}
\usepackage{amssymb}
\usepackage{amsthm}
\usepackage{mathrsfs}
\usepackage{hhline}
\usepackage{longtable}
\usepackage{amsmath}
\usepackage{array}
\usepackage{float}
\usepackage{multirow}
\usepackage{cite}
\usepackage{enumerate}
\usepackage{xcolor}

\def\ms{\medskip}
\def\nt{\noindent}

\def\rsq{\hspace*{\fill}$\Box$\medskip}

\newtheoremstyle{de}
  {10pt}          
  {10pt}  
  {\rm}  
  {}
  {\bf}  
  {. }    
  { }    
  {}     
\theoremstyle{de}

\newtheorem{example}{Example}[section]

\newtheorem{problem}{Problem}[section]
\newtheorem{question}{Question}[section]

\newtheoremstyle{theorem}
  {10pt}          
  {10pt}  
  {\it}  
  {}
  {\bf}  
  {. }    
  { }    
  {}     
\theoremstyle{theorem}
\newtheorem{theorem}{Theorem}[section]
\newtheorem{lemma}[theorem]{Lemma}

\newtheorem{corollary}[theorem]{Corollary}

\numberwithin{equation}{section}

\graphicspath{%
    {converted_graphics/}
    {/}
}

\begingroup\makeatletter\ifx\SetFigFont\undefined%
\gdef\SetFigFont#1#2#3#4#5{%
  \reset@font\fontsize{#1}{#2pt}%
  \fontfamily{#3}\fontseries{#4}\fontshape{#5}%
  \selectfont}%
\fi\endgroup%

\definecolor{vividviolet}{rgb}{0.62, 0.0, 1.0}

\begin{document}
\begin{center}
{\mathversion{bold}\Large \bf On local antimagic chromatic number of lexicographic product graphs}

\bigskip
{\large  Gee-Choon Lau$^{a,}$\footnote{Corresponding author.}, Wai Chee Shiu$^{b}$}\\

\medskip

\emph{{$^a$}Faculty of Computer \& Mathematical Sciences,}\\
\emph{Universiti Teknologi MARA (Johor Branch, Segamat Campus),}\\
\emph{85000, Malaysia.}\\
\emph{geeclau@yahoo.com}\\

\medskip
\emph{{$^b$}Department of Mathematics,}\\
\emph{The Chinese University of Hong Kong,}\\
\emph{Shatin, Hong Kong, China.}\\
\emph{wcshiu@associate.hkbu.edu.hk}\\

\end{center}

\begin{abstract} Let $G = (V,E)$ be a connected simple graph of order $p$ and size $q$. A graph $G$ is called local antimagic if $G$ admits a local antimagic labeling. A bijection $f : E \to \{1,2,\ldots,q\}$ is called a local antimagic labeling of $G$ if for any two adjacent vertices $u$ and $v$, we have $f^+(u) \ne f^+(v)$, where $f^+(u) = \sum_{e\in E(u)} f(e)$, and $E(u)$ is the set of edges incident to $u$. Thus, any local antimagic labeling induces a proper vertex coloring of $G$ if vertex $v$ is assigned the color $f^+(v)$. The local antimagic chromatic number, denoted $\chi_{la}(G)$, is the minimum number of induced colors taken over local antimagic labeling of $G$. Let $G$ and $H$ be two vertex disjoint graphs. The {\it lexicographic product} of $G$ and $H$, denoted $G[H]$, is the graph with vertex set $V(G) \times V(H)$, and $(u,u')$ is adjacent to $(v,v')$ in $G[H]$ if $(u,v)\in E(G)$ or if $u=v$ and $u'v'\in E(H)$. In this paper, we obtained sharp upper bound of $\chi_{la}(G[O_n])$ where $O_n$ is a null graph of order $n\ge 1$. Sufficient conditions for even regular bipartite and tripartite graphs $G$ to have $\chi_{la}(G)=3$ are also obtained. Consequently, we successfully determined the local antimagic chromatic number of infinitely many (connected and disconnected) regular graphs that partially support the existence of $r$-regular graph $G$ of order $p$ such that (i) $\chi_{la}(G)=\chi(G)=k$, and (ii) $\chi_{la}(G)=\chi(G)+1=k$ for each possible $r,p,k$.  \\ 


\end{abstract}


\section{Introduction}

Let $G = (V,E)$ be a connected simple graph of order $p$ and size $q$. A bijection $f : E \to \{1,2,\ldots,q\}$ is called a {\it local antimagic labeling} of $G$ if for any two adjacent vertices $u$ and $v$, we have $f^+(u) \ne f^+(v)$, where $f^+(u) = \sum_{e\in E(u)} f(e)$, and $E(u)$ is the set of edges incident to $u$. Thus, any local antimagic labeling induces a proper vertex coloring of $G$ if vertex $v$ is assigned the color $f^+(v)$. If $f$ induces $t$ distinct colors, we say $f$ is a {\it local antimagic $t$-coloring} of $G$. The {\it local antimagic chromatic number} of $G$, denoted $\chi_{la}(G)$, is the minimum number of induced colors taken over local antimagic labelings of $G$ \cite{Arumugam}.

\ms\nt Let $G$ and $H$ be two graphs. The {\it lexicographic product} of $G$ and $H$, denoted $G[H]$, is the graph with vertex set $V(G) \times V(H)$, and $(u,u')$ is adjacent to $(v,v')$ in $G[H]$ if $(u,v)\in E(G)$ or if $u=v$ and $u'v'\in E(H)$. In particular, for $n\ge 2$, $G[O_n]$ is a graph that arises by replacing every vertex of $G$, say $u$, with $n$ vertices, say $u_i$ $(1\le i\le n)$, and every edge $uv$ of $G$ with a copy of $K_{n,n}$ with edges joining vertices $u_i$ and $v_j$ for $1\le i,j\le n$. Observe that $\chi(G) = \chi(G[O_n])$ for $n\ge 2$. In~\cite{LNS, LNS-dmgt, LauShiuNg, YBY}, the local antimagic chromatic number of the join of certain graphs are determined. We are not aware of any results on local antimagic chromatic number of other graph operations. Motivated by this, in this paper, we investigated sharp upper bounds of the local antimagic chromatic number of $G[O_n]$ where $O_n$ is a null graph of order $n\ge 1$.  Sufficient conditions for even regular bipartite and tripartite graphs $G$ to have $\chi_{la}(G)=3$ are also obtained. Consequently, we successfully determined the local antimagic chromatic number of infinitely many (connected and disconnected) regular graphs that partially support the existence of $r$-regular graph $G$ of order $p$ such that (i) $\chi_{la}(G)=\chi(G)=k$, and (ii) $\chi_{la}(G)=\chi(G)+1=k$ for each possible $r,p,k$.

\ms\nt We shall need the following lemmas.

\begin{lemma}[\hspace*{-1mm}\cite{LNS-dmgt}]\label{lem-2part}
Let $G$ be a graph of size $q$. Suppose there is a local antimagic labeling of $G$ inducing a $2$-coloring of $G$ with colors $x$ and $y$, where $x<y$.  Let $X$ and $Y$ be the sets of vertices colored $x$ and $y$, respectively, then $G$ is a bipartite graph with bipartition $(X,Y)$ and $|X|>|Y|$. Moreover,
\begin{equation*} x|X|=y|Y|= \frac{q(q+1)}{2}.\end{equation*}
\end{lemma}

\nt Since bipartition of a connected bipartite graph is unique, Lemma~\ref{lem-2part} implies that

\begin{corollary}\label{cor-2part} Suppose $G$ is a connected bipartite graph of $q$ edges with bipartition $(V_1,V_2)$. If $\chi_{la}(G)=2$, then $|V_1|\ne |V_2|$ and $\binom{q+1}{2}$ is divisible by both $|V_1|$ and $|V_2|$.\end{corollary}

\nt In what follows, let $G-e$ denotes the graph $G$ with an edge $e$ deleted.

\begin{lemma}[\hspace*{-1mm}\cite{LNS-dmgt}]\label{lem-reg}  Suppose $G$ is a $d$-regular graph of size $q$. If $f$ is a local antimagic labeling of $G$, then $g = q + 1 - f$ is also a local antimagic labeling of $G$ with $c(f)= c(g)$.  Moreover, suppose $c(f)=\chi_{la}(G)$ and if $f(e)=1$ or $f(e)=q$, then $\chi_{la}(G-e)\le \chi_{la}(G)$. \end{lemma}

\section{Lexicographic Product of Regular Graphs}\label{sec-lexi}


\begin{theorem}\label{thm-lexico} Let $f$ be a local antimagic labeling of $G$ that induces $\chi_{la}(G)$ distinct vertex colors. For $n\ge 3$, $\chi(G[O_n]) \le \chi_{la}(G[O_n]) \le \chi_{la}(G)$ if for any $u,v\in V(G)$,
\begin{enumerate}[(i)]
  \item $f^+(u) = f^+(v)$ implies that $\deg(u) = \deg(v)$, and
  \item $f^+(u)\ne f^+(v)$ implies that $f^+(u)n^3 -\frac{(n^3-n)\deg(u)}{2} \ne f^+(v)n^3 -\frac{(n^3-n)\deg(v)}{2}$.
\end{enumerate}
\nt Moreover, the last equality holds if (a) $\chi(G) = \chi_{la}(G)$, or (b) $G$ is bipartite with both partite sets of same size and $\chi_{la}(G) = 3$.
\end{theorem}

\begin{proof} Let $\{u_1, u_2, \dots, u_p\}$ and $\{x_1, x_2,\dots, x_n\}$ be the vertex lists of $G$ and $O_n$, respectively, here $p=|V(G)|$.
Let $M$ be the labeling matrix corresponding to $f$ (for definition of labeling matrix, please see \cite{Shiu1998}). Let $\Omega$ be a magic square of order $n$ and $J$ be the square matrix of order $n$ whose entries are 1. Let $\Omega_i=\Omega+(i-1)n^2J$, $1\le i\le q=|E(G)|$. Note that $r_j(\Omega_i)=c_j(\Omega_i)=k+(i-1)n^3$ for each $j$, $1\le j\le n$, where $k=\frac{n^3+n}{2}$ is the magic constant of $\Omega$. Define a labeling matrix $\mathscr M$ of $G[O_n]$ by replacing each entry of $M$ by an $n\times n$ matrix as follows:
\begin{enumerate}[(1)]
\item replace $*$ by $\bigstar$, an $n \times n$ matrix whose entries are $*$;
\item replace $i$ by $\Omega_i$ if $i$ lies in the upper triangular part of $M$;
\item replace $i$ by $\Omega_i^T$, the transpose of $\Omega_i$, if $i$ lies in the lower triangular part of $M$.
\end{enumerate}

\nt Clearly all labels in $[1, qn^2]$ are filled in $\mathscr M$ and each of them appears once in the upper triangular part of $\mathscr M$. Let the corresponding labeling be $g$.

\ms\nt Suppose $(u,x)\in V(G[O_n])$, then $u=u_l$ and $x=x_j$ for some $l,j$, $1\le l\le p$ and $1\le j\le n$. Thus the row sum of $\mathscr M$ corresponding to the vertex $(u_l, x_j)$ is
\begin{align*} & \sum_{\begin{smallmatrix}\mbox{\scriptsize $i$ appears in }\\\mbox{\scriptsize the $l$-th row of $M$}\end{smallmatrix}}  \hspace*{-5mm}r_j(\Omega_i)  =
\sum_{\begin{smallmatrix}\mbox{\scriptsize $i$ appears in }\\\mbox{\scriptsize the $l$-th row of $M$}\end{smallmatrix}} \hspace*{-5mm} [k+(i-1)n^3]\\ & =\deg(u_l)k+[f^+(u_l)-\deg(u_l)]n^3 =f^+(u)n^3 -\frac{(n^3-n)\deg(u)}{2}.\end{align*}

\nt Condition (ii) implies that $g$ is a local antimagic labeling of $G[O_n]$. Condition (i) implies that $c(g)\le c(f)$. Thus $\chi_{la}(G[O_n])\le \chi_{la}(G)$.

\ms\nt Clearly, if (a) holds, then $\chi(G[O_n]) = \chi_{la}(G)$ and equality holds. Furthermore, if (b) holds, then $G[O_n]$ is also bipartite with both partite sets of same size. Lemma~\ref{lem-2part} then implies that $3\le \chi_{la}(G[O_n])\le \chi_{la}(G) = 3$. This completes the proof.
\end{proof}

\begin{example} To demonstrate the above proof, we consider $C_4[O_3]$.
Clearly the following matrices $M$ and $\Omega$  are labeling matrices corresponding to a local antimagic labeling of $C_4$ and a magic square of order 3, respectively. Note that $\chi_{la}(C_4)=3$.
\[M=\begin{pmatrix}
* & 1 & * & 4\\
1 & * & 2 & *\\
* & 2 & * & 3\\
4 & * & 3 & *
\end{pmatrix}, \quad \Omega=\begin{pmatrix}
8 & 1 & 6\\
3 & 5 & 7\\
4 & 9 & 2
\end{pmatrix}.\]
Now
\[\Omega_1=\Omega, \quad \Omega_2=\begin{pmatrix}
17 & 10 & 15\\
12 & 14 & 16\\
13 & 18 & 11
\end{pmatrix}, \quad\Omega_3=\begin{pmatrix}
26 & 19 & 24\\
21 & 23 & 25\\
22 & 27 & 20
\end{pmatrix}, \quad\Omega_4 =\begin{pmatrix}
35 & 28 & 33\\
30 & 32 & 34\\
31 & 36 & 29
\end{pmatrix}.
\]
Then \[\mathscr M=\begin{pmatrix} \bigstar & \Omega_1 & \bigstar & \Omega_4\\
\Omega_1^T & \bigstar & \Omega_2 & \bigstar\\
\bigstar & \Omega_2^T & \bigstar & \Omega_3\\
\Omega_4^T & \bigstar & \Omega_3^T & \bigstar
\end{pmatrix}=\left(\begin{array}{ccc|ccc|ccc|ccc}
* & * & * & 8 & 1 & 6 & * & * & * & 35 & 28 & 33\\
* & * & * & 3 & 5 & 7 & * & * & * & 30 & 32 & 34\\
* & * & * & 4 & 9 & 2 & * & * & * & 31 & 36 & 29\\\hline
8 & 3 & 4 & * & * & * & 17 & 10 & 15 & * & * & *\\
1 & 5 & 9 & * & * & * & 12 & 14 & 16 & * & * & *\\
6 & 7 & 2 & * & * & * & 13 & 18 & 11 & * & * & *\\\hline
* & * & * & 17 & 12 & 13 & * & * & * & 26 & 19 & 24\\
* & * & * & 10 & 14 & 18 & * & * & * & 21 & 23 & 25\\
* & * & * & 15 & 16 & 11 & * & * & * & 22 & 27 & 20\\\hline
35 & 30 & 31 & * & * & * & 26 & 21 & 22 & * & * & *\\
28 & 32 & 36 & * & * & * & 19 & 23 & 27 & * & * & *\\
33 & 34 & 29 & * & * & * & 24 & 25 & 20 & * & * & *
\end{array}
\right)\]
is a local antimagic labeling of $C_4[O_3]$. Clearly, $\chi_{la}(C_4[O_3])=3$.\rsq
\end{example}

\nt Note that if $G$ is a bipartite graph with both partite sets of same size, then $G[O_n]$ is also a bipartite graph with both partite set of same size. Since the conditions of Theorem~\ref{thm-lexico} hold for all regular graphs, we have the following corollary.

\begin{corollary}\label{cor-reglexico} Let $G$ be a regular graph. For $n\ge 3$, $\chi(G[O_n]) \le \chi_{la}(G[O_n])\le \chi_{la}(G)$. Particularly, $\chi_{la}(G[O_n]) = \chi_{la}(G)$ if (i) $\chi(G) = \chi_{la}(G)$, or (ii) $G$ is bipartite with partite sets of same size and $\chi_{la}(G) = 3$.
\end{corollary}

\nt Note that $K_t[O_n]$ is the complete $t$-partite graph with each partite set of size $n$.

\begin{corollary} For $n, t\ge 3$, $\chi_{la}(K_t[O_n]) = t$. \end{corollary}

\nt Since $\chi_{la}(K_2)$ does not exist, Corollary~\ref{cor-reglexico} does not apply to $K_2[O_n] = K_{n,n}$. Interested reader may refer to~\cite[Theorem 10]{LNS} for the exact value of $\chi_{la}(K_{n,n})$.

\ms\nt In~\cite{LauShiu}, the authors completely determined the local antimagic chromatic number of all connected regular graphs of order at most 8. Suppose $m\ge 3$. The following are regular graphs $G\ne K_t[O_n], t\ge 2, n \ne 2$, with $\chi_{la}(G)=\chi(G)$, or else $G$ is a bipartite graph with partite sets of same size with $\chi_{la}(G)=\chi(G)+1=3$. 

\begin{enumerate}[(1)]
  \item 2-regular $G=C_m$ (see~\cite{Arumugam}) with $\chi_{la}(G) = 3$.
  \item $m$-regular $G=K_{m,m}$ (see~\cite{Arumugam, LNS}) with $\chi_{la}(G)=3$ for $m\ge 2$.
  \item $m$-regular $G=C_m\vee O_{m-2}$ (see~\cite{LNS-dmgt}) with $\chi_{la}(G) = \begin{cases} 3 & \mbox{ for even $m$}, \\ 4 & \mbox{ otherwise}.\end{cases}$
  \item $(m+2)$-regular $G=C_m\vee C_m$ (see~\cite{LNS-dmgt}) with $\chi_{la}(G) = 6$ for odd $m\ge 5$.
  \item 3-regular $G=M_{2m}$ (see~\cite{LNS-dmgt}) with $\chi_{la}(G) = 3$ for odd $m$.
\end{enumerate}

\nt In~\cite[Theorem 3.15]{LauShiu}, the authors also proved the existence of non-complete regular graphs of arbitrarily large order, regularity and local antimagic chromatic numbers. By Corollary~\ref{cor-reglexico}, we immediately have the following theorem.

\begin{theorem} Let $G$ be a graph in (1) to (5) above or in~\cite[Theorem 3.15]{LauShiu}, then $\chi_{la}(G[O_n])=\chi_{la}(G)=\chi(G)$, except $\chi_{la}(G[O_n])=3=\chi_{la}(G) = \chi(G)+1$ where $G=C_m, m\ge 3$ is even, or $G=K_{m,m}, m\ge 2$, for $n\ge 3$.\end{theorem}



\nt Note that for $n\ge 3$ and the values of $m$ as in (1) to (5) above, we have
\begin{enumerate}[(1)]
  \item $C_m[O_n]$ is a $2n$-regular graph of order $mn$ with $\chi_{la}(C_m[O_n]) = \chi(C_m[O_n])+1=3$ for even $m$, or else $\chi_{la}(C_m[O_n]) = \chi(C_m[O_n])=3$. 
  \item $K_{m,m}[O_n]$ is an $mn$-regular graph of order $2mn$ with $\chi_{la}(K_{m,m}[O_n]) = \chi(K_{m,m}[O_n])+1=3$.
  \item $(C_m\vee O_{m-2})[O_n]$ is an $mn$-regular graph of order $n(2m-2)$ with $\chi_{la}((C_m\vee O_{m-2})[O_n]) = \chi((C_m\vee O_{m-2})[O_n]) = \begin{cases} 3 & \mbox{ if $m$ is even,} \\ 4 & \mbox{ otherwise.}\end{cases}$. 
  \item $(C_m \vee C_m)[O_n]$ is an $n(m+2)$-regular graph of order $2mn$ with $\chi_{la}((C_m \vee C_m)[O_n])=\chi((C_m \vee C_m)[O_n])=6$.
  \item $M_{2m}[O_n]$ is a $3n$-regular graph of order $2mn$ with $\chi_{la}(M_{2m}[O_n]) = \chi(M_{2m}[O_n])+1=3$.
\end{enumerate}

\nt Let $G$ be a non-regular graph with $\chi_{la}(G)=t\ge 2$ such that all vertices in the same color class induced by a local antimagic $t$-coloring of $G$ are of same degree.  The local antimagic chromatic number of such non-regular graphs with their induced vertex colors $X$ and vertex degree $d$ are obtained as follows.
\begin{enumerate}[(a)]
  \item $\chi_{la}(W_{4k})=\chi(W_{4k})=3$ with $(X,d)=\{(11,3),(15,3),(20,4)\,|\,k=1\}$ and $(X,d)=\{(9k+2,3),(11k+3,3),(2k(12k+1),4k)\,|\,k\ge 2\}$ (see~\cite[Theorem 5]{LNS}).
  \item $\chi_{la}(W_n)=\chi(W_m)=3$ for $n\equiv2\pmod{4}$ with $(X,d)\in\{((9m+6)/4,3),((11m+6)/4,3),(m(3m+1)/2,m)\}$ (see~\cite{Arumugam}).
  \item $\chi_{la}(W_n)=\chi(W_m)=4$ for $n\equiv1\pmod{4}$ with $(X,d)\in\{((11m+13)/4,3),((9m+11)/4,3),((5m+11)/4,3),((6m^2+m+1)/4,m)\}$ (see~\cite{Arumugam}, with minor correction).
  \item $\chi_{la}(W_n)=\chi(W_m)=4$ for $n\equiv3\pmod{4}$ with $(X,d)\in\{((9m+9)/4,3),((11m+7)/4,3),(2m+2,3),(m(3m+1)/2,m)\}$ (see~\cite{Arumugam}, with minor correction).
\end{enumerate}

\nt It is easy to check the conditions of Theorem~\ref{thm-lexico}.

\begin{corollary} For $m,n\ge 3$, $\chi_{la}(W_m[O_n])=\begin{cases}3 & \mbox{ if $m\ge 4$ is even},\\ 4 & \mbox{ otherwise.}\end{cases}$ \end{corollary}

\ms\nt Arumugam et al. in \cite{Arumugam} provided a local antimagic labelings $f$ of $C_n=x_1x_2\cdots x_nx_1$ with
\[f^+(x_i)=\begin{cases}
2n-\lfloor\frac{n}{2}\rfloor & \mbox{if $i=1$};\\
n & \mbox{if $i$ is odd, $i\ne 1$};\\
n+1 & \mbox{if $i$ is even}.\end{cases}\]
Moreover, $f(x_1x_2)=n$. We shall use this result to construct a local antimagic labeling for a $2m$-regular bipartite graph.

\begin{theorem}\label{thm-even-reg-bipartite}
Suppose $G$ is a $2m$-regular bipartite graph, $m\ge 1$, then $\chi_{la}(G)=\chi_{la}(G-e)=3$ for any edge $e$ in $G$.
\end{theorem}

\begin{proof} Since $G$ is a $2m$-regular bipartite graph, $G$ is Eulerian and both partite sets are the same size. Thus $\chi_{la}(G)\ge 3$.

\ms\nt Let $P=x_1x_2\cdots x_qx_1$ be an Euler tour of $G$, where $q$ is the size of $G$. Note that $q$ is even. We view vertices in $P$ are distinct. We apply the local antimagic labeling $f$ provided in \cite{Arumugam} on $P$ such that $f^+(x_1)=2q-\frac{q}{2}$.  We define a labeling $F:E(G)\to [1,q]$ is the same as $f$ defined on $E(P)$.

\nt Each vertex $v\in V(G)$ appears $m$ times in the trail $x_1x_2\cdots x_q$. Suppose $v=x_{i_1}=\cdots= x_{i_m}$, $1\le i_1<\cdots<i_m\le q$.  Since $G$ is bipartite, all $i_j$ are of the same parity, $1\le j\le m$.
When $P$ passes through the vertex $x_{i_j}$, it contributes a weight $f^+(x_{i_j})$ to $v$. Thus $F^+(v)=\sum_{j=1}^m f^+(x_{i_j})$. Therefore,

\[F^+(v)=\begin{cases}2q-\frac{q}{2} + (m-1)q& \mbox{if $i_1=1$};\\
m(q+1) & \mbox{if $i_1$ is even};\\
mq & \mbox{if $i_1$ is odd, $i_1\ne 1$}.\\
\end{cases}\]
So $F$ is a local antimagic $3$-coloring for $G$ and hence $\chi_{la}(G)=3$.

\ms \nt Now, for any edge $e$ in $G$, we may choose an Euler tour $P=x_1x_2\cdots x_q$ such that $e=x_1x_2$.  Now $F(x_1x_2)=q$. By Lemma~\ref{lem-reg}, we have $\chi_{la}(G-e)\le 3$ and hence $\chi_{la}(G-e)=3$.
\end{proof}

\begin{example} For $m,n\ge 2$, let $G_{m,n}$ be a graph with vertex set \[V(G_{m,n}) = \{v_{i,j}\,|\, 1\le i\le m, 0\le j\le 2n-1\}\] and edge set \[E(G_{m,n}) = \{v_{1,j}v_{1,j+1}, v_{i,j}v_{i+1,j+1}, v_{i+1,j}v_{i,j+1}, v_{m,j}v_{m,j+1} \,|\, 1\le i\le m-1, 0\le j\le 2n-1\},\]
here $v_{i,2n}=v_{i,0}$ by convention. Note that $G_{m,n}$ is a $4$-regular bipartite graph. By Theorem~\ref{thm-even-reg-bipartite}, $\chi_{la}(G_{m,n})=\chi_{la}(G_{m,n}-e)=3$ for any edge $e$ in $G_{m,n}$.\rsq
\end{example}

\begin{corollary}\label{cor-product-mcycles} For $m\ge1$ and $n_i\ge 1$, $1\le i\le m$, let $G=C_{2n_1}\times \cdots \times C_{2n_m}$. Suppose $n\ge 2$, then $\chi_{la}(G)=\chi_{la}(G[O_n])=3$. 
\end{corollary}

\nt Observe that $C_{2n}[O_2]$ is a 4-regular bipartite graph. By Theorem~\ref{thm-even-reg-bipartite}, we have

\begin{corollary} For $n\ge 2$, $\chi_{la}(C_{2n}[O_2])=3$. \end{corollary}

\nt Note that we can also extend the argument for Theorem~\ref{thm-even-reg-bipartite} to certain tripartite graphs $G$ with $\chi(G)=3$. Let us consider an Eulerian tripartite graph $G$ of size $q$ as follows:

\begin{enumerate}[(T1)]
  \item The 3 partite sets are $V_1=\{w\}$, $V_2=\{u_i\,|\,1\le i\le x, x\ge 2\}$, $V_3=\{v_j\,|\,1\le j\le y, y\ge 2\}$. 
  \item $w$ is adjacent to $2a$ vertices in $V_2$ and $2b$ vertices in $V_3$ with $a,b \ge 1$, each vertex in $V_2$ has degree $2m$, each vertex in $V_3$ has degree $2n$. 
\end{enumerate}

\nt By counting all edges incident to vertices in $V_2$, and all edges incident to $w$ and vertices in $V_3$, we get $q=2mx+2b$. Similarly, we have $q=2ny+2a$. Thus $q$ is even and $mx+b=ny+a$.

\ms \nt Suppose $C$ is an Eulerian tour of $G$, then $C$ is a union of some edge-disjoint closed trails. We may assume that each such trail contains $w$ once beginning and ending at $w$. Let $R_i$, $S_j$, $T_k$ be such closed trails with the following property:
\begin{enumerate}[1.]
\item $C=(\bigcup\limits_{1=1}^\alpha R_i)(\bigcup\limits_{j=1}^\beta S_j)(\bigcup\limits_{k=1}^\gamma T_k)$, where $\alpha, \beta, \gamma\ge 0$. By convention, empty union does not appear in the formula.
\item In $R_i$, $w$ is adjacent to two vertices in $V_2$.
\item In $S_j$, $w$ is adjacent to two vertices in $V_3$.
\item In $T_k$, $w$ is adjacent to one vertex in $V_2$ and one vertex in $V_3$.
\end{enumerate}
Note that $R_i$ and $S_j$ are of even lengths, $T_k$ is of odd length. Since $G$ is not bipartite, $G$ contains at least one odd cycle. Since $q$ is even, $\gamma$ is a positive even number. Moreover, $2\alpha+\gamma=2a$ and $2\beta+\gamma=2b$.

\ms\nt We assume the second vertex of $T_k$ is in $V_2$ and the last second vertex is in $V_3$ when $k$ is odd; whereas the second vertex of $T_k$ is in $V_3$ and the last second vertex is in $V_2$ when $k$ is even. We can now rearrange the sequence of visit of the Eulerian tour as follows:
\begin{equation}\label{eq-eulerian-even} R_1\cdots R_\alpha T_1\cdots T_{\gamma-1} S_1\cdots S_\beta T_\gamma.\end{equation}
Thus, the Eulerian tour $C$ have the property that
\begin{equation}\label{cond-C} \mbox{the two vertices immediately before and after $w$ in the tour are in the same partite set.}\tag{\dag}\end{equation}

\begin{theorem}\label{thm-tripartite-even} Suppose $G$ is the tripartite graph satisfying Conditions (T1) and (T2). Then $\chi_{la}(G) = 3$. \end{theorem}

\begin{proof} Let $C=x_1\cdots x_{q}x_1$ be the Eulerian tour and be arranged as in~\eqref{eq-eulerian-even}, i.e., $C$ satisfies the Condition~\eqref{cond-C}. We shall say that the vertex $x_j$ is an odd (resp. even) term of $C$ if $j$ is odd (resp. even). By~\eqref{eq-eulerian-even} and Condition~\eqref{cond-C}, all vertices in $V_2$ are even terms and all vertices in $V_3$ are odd terms. Therefore, vertex $w$ appears as an odd term for $\alpha+\gamma/2=a$ times and as an even term for $\beta+\gamma/2=b$ times.




\ms\nt Let $f$ be the edge labeling defined in~\cite{Arumugam} for a cycle as stated before Theorem~\ref{thm-even-reg-bipartite} above. We have $f^+(u_i) = m(q+1)$ for each $u_i\in V_2$, $f^+(v_j) = nq$ for each $v_j\in V_3$ and $f^+(w) = (a-1)q + b(q+1) +3q/2 = (a+b+\frac{1}{2})q +b=(a+b)q + b + ny + a=(a+b)(q+1)+ny$.

\ms\nt We now have the followings:
\begin{enumerate}[(a)]
  \item Since $q> m$, $f^+(u_i)-f^+(v_j)=(m-n)q+m\ne 0$. Thus $f^+(u_i) \ne f^+(v_j)$.
  \item $f^+(v_j)-f^+(w)=(n-a-b-\frac{1}{2})q-b$. Since $q>2b$, $f^+(v_j)-f^+(w)\ne 0$. Thus $f^+(v_j)\ne f^+(w)$.
    \item $f^+(u_j)-f^+(w)=(m-a-b)(q+1)-ny$. If $f^+(u_j)=f^+(w)$, then $(m-a-b)(q+1)=ny=\frac{q}{2}-a>0$. Hence the left hand side of this equation is greater than $q$. But the right hand side of this equation is less than $q/2$ which is a contradiction.
\end{enumerate}

\nt Thus, $G$ admits a local antimagic labeling that induces 3 distinct vertex labels. Since $\chi_{la}(G)\ge \chi(G) = 3$, the theorem holds.
\end{proof}




\nt Let us consider the case that $w$ is adjacent to $2a+1$ vertices in $V_2$ and $2b+1$ vertices in $V_3$. That means we change Condition~(T2) as
\begin{enumerate}[(T1)]
\item[(T2$'$)] $w$ is adjacent to $2a+1$ vertices in $V_2$ and $2b+1$ vertices in $V_3$ with $a,b \ge 0$, each vertex in $V_2$ has degree $2m$, each vertex in $V_3$ has degree $2n$.
\end{enumerate}

\begin{theorem}\label{thm-tripartite-odd} Suppose $G$ is the tripartite graph satisfying Conditions (T1) and (T2$'$). Then $\chi_{la}(G) = 3$. \end{theorem}
\begin{proof} Keep all notation used in the proof of Theorem~\ref{thm-tripartite-even}. Now $q=2mx+2b+1=2ny+2a+1$, and hence $\gamma$ is a positive odd number. Moreover, $2\alpha+\gamma=2a+1$ and $2\beta+\gamma=2b+1$. The sequence of visit of an Eulerian tour $C$, which begins and ends at $w$ is
\begin{equation}\label{eq-eulerian-odd} R_1\cdots R_\alpha T_1\cdots T_{\gamma} S_1\cdots S_\beta.\end{equation}
All intermediate vertices of $C$ satisfy Condition~\eqref{cond-C}. Thus all vertices in $V_2$ are even terms and all vertices in $V_3$ are odd terms. Therefore, if we do not count the first and the last $w$ in $C$, $w$ appears as an odd term for $\alpha+(\gamma-1)/2=a$ times and as an even term for $\beta+(\gamma-1)/2=b$ times.



\ms\nt Let $f$ be the edge labeling defined in~\cite{Arumugam} for a cycle as stated before Theorem~\ref{thm-even-reg-bipartite} above. Then
$f^+(u_i) = m(q+1)$ for each $u_i\in V_2$, $f^+(v_j) = nq$ for each $v_j\in V_3$ and $f^+(w) = aq + b(q+1) +2q-\lfloor q/2\rfloor = (a+b+3/2)q + b + 1/2 = (a+b+1)(q+1)+ny$.

\ms\nt We now have the followings:
\begin{enumerate}[(a)]
\item By the same argument in the proof of Theorem~\ref{thm-tripartite-even} we have $f^+(u_i)\ne f^+(v_j)$.
\item $f^+(v_j)-f^+(w)=(n-a-b-\frac{3}{2})q-b-\frac{1}{2}$. If $f^+(v_j)=f^+(w)$, then $(n-a-b-\frac{3}{2})q=b+\frac{1}{2}\ne 0$. Thus $|n-a-b-\frac{3}{2}|\ge \frac{1}{2}$. Hence $|n-a-b-\frac{3}{2}|q>\frac{1}{2}(2b+1)=b+\frac{1}{2}$ which is a contradiction. Therefore, $f^+(v_j)\ne f^+(w)$.
\item $f^+(u_i)-f^+(w)=(m-a-b-1)(q+1)-ny$. If $f^+(u_i)=f^+(w)$, then $(m-a-b-1)(q+1)=ny$. Since $ny>0$, $(q+1)\le (m-a-b-1)(q+1)=ny\le \frac{q-1}{2}$ which is a contradiction. Hence $f^+(u_i)\ne f^+(w)$.
\end{enumerate}
\nt Thus, $G$ admits a local antimagic labeling that induces 3 distinct vertex labels. Since $\chi_{la}(G)\ge \chi(G) = 3$, the theorem holds.
\end{proof}

\begin{theorem}\label{thm-lexico-tripartite-even} Let $G$ be the tripartite graph as defined in Theorem~\ref{thm-tripartite-even}. If $t\ge 3$, then $\chi_{la}(G[O_t]) = 3$. \end{theorem}

\begin{proof} Suppose $u$ and $v$ are two distinct vertices of $G$. Let $f$ be the local antimagic labeling of $G$ as defined in the proof of Theorem~\ref{thm-tripartite-even}. Clearly, $f^+(u) = f^+(v)$ implies that $deg(u) = deg(v)$. Consider the vertices $w$, $u_i$ and $v_j$ of $G$ with $f^+(w) = (a+b)(q+1)+ny$, $f^+(u_i)=m(q+1)$, and $f^+(v_j) = nq$, with $\deg(w) = 2a+2b$, $\deg(u_i) = 2m$ and $\deg(v_j)=2n$.

\ms\nt Consider $w$ and $u_i$. Suppose $f^+(w)t^3 -\frac{(t^3-t)\deg(w)}{2} = f^+(u_i)t^3 -\frac{(t^3-t)\deg(u_i)}{2}$, then $$[(a+b)(q+1)+ny]t^2 - (t^2-1)(a+b) = m(q+1)t^2 - (t^2-1)m.$$ This means $[(a+b-m)(q+1)+ny]t^2 = (t^2-1)(a+b-m)$ so that $(m-a-b)(qt^2+1) = nyt^2>0$.
Hence $m>a+b$. Since $q>2ny$, $(m-a-b)(qt^2+1)>nyt^2$, a contradiction.

\ms\nt Consider $w$ and $v_j$. Suppose $f^+(w)t^3 -\frac{(t^3-t)\deg(w)}{2} = f^+(v_j)t^3 -\frac{(t^3-t)\deg(v_j)}{2}$, then $$[(a+b)(q+1)+ny]t^2 - (t^2-1)(a+b) = nqt^2 - (t^2-1)n.$$ This means $[(a+b-n)q+(a+b)+ny]t^2 = (t^2-1)(a+b-n)$ so that  $(n-a-b)(qt^2+1)=n(y+1)t^2>0$. Thus $n>a+b$. Since $q>2ny\ge ny+n$,  $(n-a-b)(qt^2+1)>n(y+1)t^2$, a contradiction.

\ms\nt Consider $u_i$ and $v_j$. Suppose $f^+(u_i)t^3 -\frac{(t^3-t)\deg(u_i)}{2} = f^+(v_j)t^3 -\frac{(t^3-t)\deg(v_j)}{2}$, then $$m(q+1)t^2 - (t^2-1)m = nqt^2 - (t^2-1)n.$$ This means  $(n-m)(qt^2+1)=nt^2>0$. Thus $n>m$. Since $q>2n$,  $(n-m)(qt^2+1)>nt^2$ a contradiction.

\ms\nt Since $\chi(G)=\chi_{la}(G)=3$, by Theorem~\ref{thm-tripartite-even}, we have $\chi_{la}(G[O_t]) = 3$ for $t\ge 3$.
\end{proof}

\begin{theorem}\label{thm-lexico-tripartite-odd} Let $G$ be the tripartite graph as defined in Theorem~\ref{thm-tripartite-odd}. If $t\ge 3$, then $\chi_{la}(G[O_t]) = 3$. \end{theorem}

\begin{proof} Suppose $u$ and $v$ are two distinct vertices of $G$. Let $f$ be the local antimagic labeling of $G$ as defined in the proof of Theorem~\ref{thm-tripartite-odd}. Clearly, $f^+(u) = f^+(v)$ implies that $deg(u) = deg(v)$. Consider the vertices $w$, $u_i$ and $v_j$ of $G$ with $f^+(w) = (a+b+1)(q+1)+ny$, $f^+(u_i)=m(q+1)$, and $f^+(v_j) = nq$, with $\deg(w) = 2a+2b+2$, $\deg(u_i) = 2m$ and $\deg(v_j)=2n$. Similar to the argument in the proof of Theorem~\ref{thm-lexico-tripartite-even}, we can conclude that $\chi_{la}(G[O_t]) = 3$ for $t\ge 3$.
\end{proof}


\section{Disconnected regular graphs}

In~\cite{Baca}, the authors obtained the following two results on local antimagic chromatic number of disconnected regular graphs.

\begin{theorem}\label{thm-4r} Let $G$ be a $4r$-regular graph, $r \ge 1$. For every positive integer $m$,
$\chi_{la}(mG) \le \chi_{la}(G)$. \end{theorem}

\begin{theorem}\label{thm-4r+2} Let $G$ be a $(4r+2)$-regular graph, $r\ge 0$, containing a $2$-factor consisting only from even cycles. For every positive integer $m$, $\chi_{la}(mG)\le \chi_{la}(G)$. \end{theorem}

\begin{corollary}  Let $n\ge 2, m\ge 1$, then $\chi_{la}(mC_{2n}) = 3$.\end{corollary}

\nt By Corollary~\ref{cor-reglexico}, we immediately have

\begin{corollary} Let $n\ge 2, m,k\ge 1$, $\chi_{la}(mC_{2n}[O_k])=3$. \end{corollary}

\nt Recalled that the local antimagic chromatic number of the following even regular graphs $G$ are known:

\begin{enumerate}
  \item $G = K_{2n+1}$ with $\chi(G)=\chi_{la}(G)=n$ for $n\ge 2$ (see~\cite{Arumugam}).
  \item $G = K_{n,n}$ with $\chi(G)=2, \chi_{la}(G) = 3$ for $n\ge 2$ (see~\cite{Arumugam, LNS}).
  \item $G = C_n\vee O_{n-2}$ with $\chi_{la}(G) = 3$ for even $n\ge 4$ (see~\cite{LNS-dmgt}).
  \item $G$ is an even regular graph of order at most 8 with $\chi(G) = \chi_{la}(G)$ as in~\cite[Section 2]{LauShiu}.
  \item $G = G_k, k\ge 1$, is an even regular graph of order $p\ge 5$ with $\chi(G) = \chi_{la}(G)$ as in~\cite[Corollary 3.6]{LauShiu}.
  \item $G = G_k, k\ge 1$, is an even regular graph of order $p\ge 3$ with $\chi(G) = \chi_{la}(G)$ as in~\cite[Corollary 3.8]{LauShiu}.
  \item $G$ is an even regular graph with $\chi(G) = \chi_{la}(G)$ as in~\cite[Theorem 3.15]{LauShiu}.
\end{enumerate}

\begin{theorem} For every positive integer $m$ and $n\ge 3$, if $G$ is an even regular graph in (1) to (7) above satisfying Theorems~\ref{thm-4r} and~\ref{thm-4r+2}, then $\chi_{la}(mG) = \chi_{la}(mG[O_n])= \chi_{la}(G) = \chi(G)$.
\end{theorem}

\section{Conclusion and Open problems}

In this paper, we obtained sufficient conditions for $\chi_{la}(G[O_n]) \le \chi_{la}(G)$. Sufficient conditions for even regular bipartite and tripartite graphs $G$ to have $\chi_{la}(G)=3$ are also obtained. Consequently, we successfully determined the local antimagic chromatic number of infinitely many (connected and disconnected) regular graphs. Note that $$\chi_{la}(C_m[O_n]) = \begin{cases} \chi_{la}(C_m)=\chi(C_m)+1=3 & \mbox{ for even $m\ge 4$ and $n\ge 2$,}\\ \chi_{la}(C_m) = \chi(C_m) = 3 & \mbox{ for odd $m\ge 3$ and $n\ge 2$,}\end{cases}$$
and that $C_m[O_n]$ is $2n$-regular of order $mn$ for $m\ge 3, n\ge 1$. Thus, together with those results in~\cite{LauShiu}, we partially answer the following question.

\begin{question} Does there exist $r$-regular graph $G$ of order $p$ such that (i) $\chi_{la}(G)=\chi(G)=k$, and (ii) $\chi_{la}(G)=\chi(G)+1=k$ for each possible $r,p,k$?\end{question}

\ms\nt The following problems arise naturally.

\begin{problem} For any connected graph $G$, obtain a sharp upperbound for $\chi_{la}(G[O_2])$. \end{problem}

\begin{problem} For any connected graphs $G$ and $H\not= O_n, n\ge 2$, obtain a sharp upperbound for $\chi_{la}(G[H])$. \end{problem}

\begin{problem} For any connected graphs $G$ and $H$, investigate $\chi_{la}(G\times H)$. \end{problem}

%
%
%


\end{document}